\documentclass[review]{elsarticle}
\makeatletter
\def\ps@pprintTitle{%
	\let\@oddhead\@empty
	\let\@evenhead\@empty
	\let\@oddfoot\@empty
	\let\@evenfoot\@oddfoot
}
\makeatother

\usepackage{lineno}
\usepackage{hyperref}
\usepackage[dvipsnames]{xcolor}
\newcommand\myshade{85}
\colorlet{mylinkcolor}{blue}
\colorlet{mycitecolor}{red}
\colorlet{myurlcolor}{blue}
\hypersetup{
	linkcolor  = mylinkcolor!\myshade!black,
	citecolor  = mycitecolor!\myshade!black,
	urlcolor   = myurlcolor!\myshade!black,
	colorlinks = true,
}
\usepackage{yuyuan, braket, geometry}
\usepackage{mathtools}
\mathtoolsset{showonlyrefs}
\usepackage{xspace}

\modulolinenumbers[5]

\def\turan{Tur\'an\xspace}
\newcommand{\pk}[1]{P^{(\alpha_{#1},\beta_{#1})}}

\newcommand{\dx}{\dfrac{d}{dx}}
\newcommand{\np}{_{n-1}}
\newcommand{\nn}{_{n+1}}

\newcommand{\drv}[1]{\left(#1\right)^\prime}
\newcommand{\pkk}{P_n^n}
\newcommand{\pkn}{P_n^{n+1}}
\newcommand{\ppk}{P_{n-1}^n}
\newcommand{\pnn}{P_{n+1}^{n+1}}

\begin{document}

\begin{frontmatter}


\title{An inequality for Jacobi polynomials of form $P_n^{(\alpha_n,\beta_n)}(x)$}


%
%

\author{Zhulin He\corref{mycorrespondingauthor}}
\address{Department of Statistics, Iowa State University, Ames, IA 50011, USA}
\ead{hezhulin@iastate.edu}
\cortext[mycorrespondingauthor]{Corresponding author}

\author{Yuyuan Ouyang\corref{}}
\address{Department of Mathematical Sciences, Clemson University, Clemson, SC 29634, USA}

\begin{abstract}
We prove an inequality for Jacobi polynomials that 
\begin{align}
\Delta_n(x):=P_n^{(\alpha_n,\beta_n)}(x)P_n^{(\alpha\nn,\beta\nn)}(x)- P\np^{(\alpha_n,\beta_n)}(x)P\nn^{(\alpha\nn,\beta\nn)}(x)\le 0,\ \forall x\ge 1,
\end{align}
where $\alpha_n=an$ and $\beta_n=bn$ for some $a,b\ge 0$. The above inequality has a similar taste as the Tu\'ran type inequalities, but with $\alpha_n$ and $\beta_n$ that depends linearly on $n$.
\end{abstract}

\begin{keyword}
\turan inequality\sep Jacobi polynomials
\end{keyword}

\end{frontmatter}

\linenumbers

\section{Introduction}
The Jacobi polynomials $P_n^{(\alpha,\beta)}(x)$ are a class of orthogonal polynomials that are well studied in many literatures. The polynomial representation with real variable $x$ is
\begin{align}
\label{eq:Jacobi}
P_n^{(\alpha,\beta)}(x) = \sum_{t=0}^{n}\binom{n+\alpha}{n-t}\binom{n+\beta}{t}\left(\frac{x-1}{2}\right)^t\left(\frac{x+1}{2}\right)^{n-t}.
\end{align}
In \cite{gasper1972inequality}, it was proved that
\begin{align}
\label{eq:GasperResult}
\Delta_n(x):=R_n^2(x) - R_{n-1}(x)R_{n+1}(x)\ge \frac{(\beta-\alpha)(1-x)}{2(n+\alpha+1)(n+\beta)}R_n^2(x),\ \forall x\in[-1,1], n\ge 1,\alpha,\beta>-1,
\end{align}
where
\begin{align}
R_n(x):=\frac{P_n^{(\alpha,\beta)}(x)}{P_n^{(\alpha,\beta)}(1)}.
\end{align}
Consequently, we have 
$\Delta_n(x)\ge 0$ 
for all $x\in[-1,1], n\ge 1, \beta\ge \alpha>-1.$ Such result is known as a Tu\'ran type inequality that originates from the studies of Legendre polynomials by Tu\'ran in \cite{turan1950zeros} (see also \cite{szego1948inequality}). It should be noted that the discussion is restricted to $x\in[-1,1]$ here because that the Jacobi polynomials are orthogonal in $[-1,1]$. The definition we use in \eqref{eq:Jacobi} is in fact well defined for any $x\in\R$.

In this study, we will prove the following inequality for Jacobi polynomials:
\begin{align}
\label{eq:Turan}
\Delta_n(x):=P_n^{(\alpha_n,\beta_n)}(x)P_n^{(\alpha\nn,\beta\nn)}(x)- P\np^{(\alpha_n,\beta_n)}(x)P\nn^{(\alpha\nn,\beta\nn)}(x)\le 0,\ \forall x\ge 1.
\end{align}
Here, $\alpha_n$ and $\beta_n$ are dependent on $n$ with
\begin{align}
\label{eq:alphabetax}
\alpha_n=an, \beta_n=bn
\end{align}
for some $a,b\ge 0$.
The inequality \eqref{eq:Turan} is different from \eqref{eq:GasperResult} due to such dependence on $n$. It should also be noted that unlike \eqref{eq:GasperResult} we are not considering $x\in[-1,1]$ in \eqref{eq:Turan}. This is because that polynomials $\{P_n^{(\alpha_n,\beta_n)}(x)\}_{n=0}^{\infty}$ are in general not orthogonal on $[-1,1]$.

\section{Notations and preliminaries}
For easy reading, we will use the following notations:
\begin{align}
P_l^{k}:=\pk{k}_l(x), \drv{P_l^{k}}:=\dx\pk{k}_l(x).
\end{align}
Under the above notation, we have 
\begin{align}
\label{eq:Delta}
\Delta_n(x) = \pkk\pkn - \ppk\pnn.
\end{align} 
We will also use notations $\Delta_n$ and $\Delta_n'$ without $x$ for convenience.

The recurrence formula for differentiation (see, e.g., Section 4.5 in \cite{szeg1939orthogonal}) is an important tool for our analysis. In particular, we have
\begin{align}
(1-x^2)\drv{\pkk} = &A_n\pkk + B_n\ppk
\\
(1-x^2)\drv{\pkn} = &C_n\pkn + D_n\pnn
\end{align}
where
\begin{align}
A_n = & -n\left(x+\dfrac{\beta_n-\alpha_n}{2n+\alpha_n+\beta_n}\right)
,
\\
B_n=&\frac{2(n+\alpha_n)(n+\beta_n)}{2n+\alpha_n+\beta_n} 
,
\\
C_n=&(n+\alpha\nn+\beta\nn+1)\left(x - \frac{\beta\nn-\alpha\nn}{2n+\alpha\nn+\beta\nn+2}\right) 
,
\\
D_n = & -\frac{2(n+1)(n+\alpha\nn+\beta\nn+1)}{2n+\alpha\nn+\beta\nn+2}
.
\end{align}
Substituting the definition of $\alpha_n$ and $\beta_n$ in \eqref{eq:alphabetax} to the above equations, we have
\begin{align}
\label{eq:DiffRecur_kk_minus}
(1-x^2)\drv{\pkk} = &nA\pkk + nB\ppk
\\
\label{eq:DiffRecur_kn_plus}
(1-x^2)\drv{\pkn} = &(n+1)C\pkn + (n+1)D\pnn
\end{align}
where
\begin{align}
\label{eq:ABCD}
A =  -x - \frac{b-a}{2+a+b},
B= \frac{2(1+a)(1+b)}{2+a+b},
C=(1+a+b)\left(x - \frac{b-a}{2+a+b}\right),
\text{ and }
D =  -\frac{2(1+a+b)}{2+a+b}.
\end{align}
Note that \eqref{eq:DiffRecur_kk_minus} and \eqref{eq:DiffRecur_kn_plus} are polynomial equalities. Therefore, although most studies of Jacobi polynomials focus on the case when $x\in[-1,1]$, the relations \eqref{eq:DiffRecur_kk_minus} and \eqref{eq:DiffRecur_kn_plus} hold for all $x\in\R$.

%

\vgap

Since Jacobi polynomials $P_n^{\alpha,\beta}(x)$ with fixed $\alpha$ and $\beta$ are orthogonal polynomials, they satisfy the following important inequality (see, e.g., (3.3.6) in \cite{szeg1939orthogonal}):
\begin{align}
	(P_{n+1}^{\alpha,\beta}(x))'P_{n}^{\alpha,\beta}(x) - (P_{n}^{\alpha,\beta}(x))'P_{n+1}^{\alpha,\beta}(x) >0,\ \forall x\in\R.
\end{align}
Setting $\alpha=\alpha_{n+1}$ and $\beta = \beta_{n+1}$ in the above equation, we have
\begin{align}
\label{eq:Wronskian}
\pkn\drv{\pnn}-\pnn\drv{\pkn}> 0.
\end{align}
The above inequality will be useful in our proof of \eqref{eq:Turan}.


%

\vgap

\section{Main results}

We will prove \eqref{eq:Turan} in this section. The general idea of our proof is to show that there exists $r$ and $s$ such that at all critical points of the function $f^{r,s}_n(x):=(x-1)^r(x+1)^s\Delta_n(x)$ in $(1,\infty)$, the values of $f_n^{r,s}(x)$ are all non-positive. As a consequence, we have $f_n^{r,s}(x)\le 0$ in $[1,\infty)$, and hence $\Delta_n\le 0$. Similar idea was used in the proofs of orthogonal polynomial inequality relations in \cite{bustoz1979inequalities,pyung2004inequalities}. Noting that
\begin{align}
\label{eq:df}
\dx f_n^{r,s}(x) = & (x-1)^{r-1}(x+1)^{s-1}\left\{[r(x+1)+s(x-1)]\Delta_n + (x^2-1)\Delta_n' \right\},
\end{align}
any critical points $x>1$ can be characterized by the following relationship between $\Delta_n$ and $\Delta_n'$:
\begin{align}
[r(x+1)+s(x-1)]\Delta_n + (x^2-1)\Delta_n' =0.
\end{align}
To start with, we prove a technical lemma below that is related to the above equation.

\vgap

\begin{lem}
	\label{lem:dDelta2Delta}
	\begin{align}
	\label{eq:dDelta2Delta}
	E_n\Delta_n + (x^2-1)\Delta_n'=(x^2-1)\left[ \frac{\pkk\drv{\pkn}}{n+1} - \frac{\drv{\pkk}\pkn}{n}\right],
	\end{align}
	where
	\begin{align}
	\label{eq:E}
	E_n:=(n+1)A + nC.
	\end{align}
	\begin{proof}
		By \eqref{eq:Delta}, \eqref{eq:DiffRecur_kk_minus} and \eqref{eq:DiffRecur_kn_plus} we have
		\begin{align}
		& (1-x^2)\left[\Delta_n' -  \frac{\pkk\drv{\pkn}}{n+1} + \frac{\drv{\pkk}\pkn}{n}\right]
		\\
		= & (1-x^2)\left[\frac{n+1}{n}\drv{\pkk}\pkn+ \frac{n}{n+1}\pkk\drv{\pkn} - \drv{\ppk}\pnn - \ppk\drv{\pnn}\right]
		\\
		=& \frac{n+1}{n}(nA\pkk+nB\ppk)\pkn + \frac{n}{n+1}((n+1)C\pkn+(n+1)D\pnn)\pkk 
		\\
		&- (nC\ppk+nD\pkk)\pnn - ((n+1)A\pnn+(n+1)B\pkn)\ppk
		\\
		= & \left((n+1)A+nC\right)\pkk\pkn - \left((n+1)A + nC\right)\ppk\pnn.
		\end{align}
		We conclude \eqref{eq:dDelta2Delta} immediately from the above result and \eqref{eq:E}.
	\end{proof}
\end{lem}

\vgap

We make one observation from the above lemma. Applying \eqref{eq:ABCD} to \eqref{eq:E}, we have
\begin{align}
E_n = (an+bn-1)x + \frac{[(2+a+b)n+1](b-a)}{2+a+b}.
\end{align}
Therefore, we can obtain
\begin{align}
\label{eq:rsE}
r(x+1) + s(x-1) = E_n
\end{align}
by setting 
\begin{align}
\label{eq:rs}
\begin{aligned}
r = & \frac{1}{2}\left[(an+bn-1)x + \frac{[(2+a+b)n+1](b-a)}{2+a+b}\right],
\\
s = & \frac{1}{2}\left[(an+bn-1)x - \frac{[(2+a+b)n+1](b-a)}{2+a+b}\right].
\end{aligned}
\end{align}

\vgap

We are now ready to prove \eqref{eq:Turan}.

\vgap

\begin{thm}
	For all $x\ge 1$, we always have
	\begin{align}
	\Delta_n(x)\le 0.
	\end{align}
	Here the inequality becomes equality if and only if $x=1$.
	\begin{proof}
		By \eqref{eq:DiffRecur_kk_minus} and \eqref{eq:DiffRecur_kn_plus} we have
		\begin{align}
		& (1-x^2)\left(\dfrac{\pkk\drv{\pnn}}{n+1}-\dfrac{\pnn\drv{\pkk}}{n}\right)
		\\
		= & \pkk(A\pnn+B\pkn) - \pnn(A\pkk+B\ppk)
		\\
		=& B(\pkk\pkn - \ppk\pnn).
		\end{align}
		Using the above relation and noting the definition of $\Delta_n$ in \eqref{eq:Delta}, we have
		\begin{align}
		\label{eq:Delta2drv}
		\Delta_n = \frac{1-x^2}{B}\left(\dfrac{\pkk\drv{\pnn}}{n+1}-\dfrac{\pnn\drv{\pkk}}{n}\right).
		\end{align}
		
		Let us define
		\begin{align}
		\label{eq:f}
		f_n^{r,s}(x):=(x-1)^r(x+1)^s\Delta_n(x),
		\end{align}
		where $r$ and $s$ are defined in \eqref{eq:rs}, so that \eqref{eq:rsE} holds. Clearly, the sign of $\Delta_n(x)$ and $f_n^{r,s}(x)$ are the same in $[1,\infty)$. 
		Noting the derivative of $f_n^{r,s}$ in \eqref{eq:df}, the relations \eqref{eq:dDelta2Delta}, and \eqref{eq:rsE}, we have		
		\begin{align}
		\dx f_n^{r,s}(x) 
		= & (x-1)^{r-1}(x+1)^{s-1}\left(E_n\Delta_n - (1-x^2)\Delta_n'\right)
		\\
		= & (x-1)^{r}(x+1)^{s}\left[\frac{\pkk\drv{\pkn}}{n+1} - \frac{\drv{\pkk}\pkn}{n}\right].
		\end{align}
		Therefore, for any critical points of $f_n^{r,s}(x)$ in $(1, \infty)$, we always have
		\begin{align}
		\label{eq:dPPrelation}
		\frac{\pkk\drv{\pkn}}{n+1} = \frac{\drv{\pkk}\pkn}{n}.
		\end{align}
		
		Suppose that $x>1$ is a critical point $f_n^{a,b}(x)$ in $(1, \infty)$. Noting the definition of $f_n^{a,b}$ in \eqref{eq:f}, using the above relation, \eqref{eq:Wronskian}, and \eqref{eq:Delta2drv}, we obtain 
		\begin{align}
		f_n^{r,s}(x) = & -\frac{(x-1)^{r+1}(x+1)^{s+1}}{B}\left(\dfrac{\pkk\drv{\pnn}}{n+1}-\dfrac{\pnn\drv{\pkk}}{n}\right)
		\\
		= & -\frac{(x-1)^{r+1}(x+1)^{s+1}}{B\pkn}\left(\dfrac{\pkk\drv{\pnn}\pkn}{n+1}-\dfrac{\pnn\drv{\pkk}\pkn}{n}\right)
		\\
		= & -\frac{(x-1)^{r+1}(x+1)^{s+1}\pkk}{B(n+1)\pkn}\left(\pkn\drv{\pnn}-\pnn\drv{\pkn}\right).
		\end{align}
		Since $x>1$, by \eqref{eq:Jacobi} we have in the above relation that $\pnn,\pkn> 0$. Therefore, applying \eqref{eq:Wronskian} to the above, we obtain $f_n^{r,s}(x)< 0$ at any critical points in $(1,\infty)$. If $f_n^{r,s}(1)\le 0$, then we conclude that $f_n^{r,s}(x)\le 0$ for all $x>1$, and so $\Delta_n<0$. To finish the proof, it suffices to check the value of $f_n^{r,s}$ when $x=1$ and $x\to\infty$. In fact, by \eqref{eq:Jacobi} we have 
		\begin{align}
		\Delta_n(1) = & P_n^{\alpha_n,\beta_n}(1)P_n^{\alpha\nn,\beta\nn}(1) - P\np^{\alpha_n,\beta_n}(1)P\nn^{\alpha\nn,\beta\nn}(1)
		\\
		= & \begin{pmatrix}
		n + \alpha_n
		\\
		n
		\end{pmatrix}\begin{pmatrix}
		n + \alpha\nn
		\\
		n
		\end{pmatrix} - \begin{pmatrix}
		n - 1 + \alpha_n
		\\
		n - 1
		\end{pmatrix}\begin{pmatrix}
		n+1 + \alpha\nn
		\\
		n+1
		\end{pmatrix}
		\\
		= & \frac{n+\alpha_n}{n}\begin{pmatrix}
		n + \alpha_n - 1
		\\
		n - 1
		\end{pmatrix}\begin{pmatrix}
		n + \alpha\nn
		\\
		n
		\end{pmatrix} - \dfrac{n+1+\alpha\nn}{n+1}\begin{pmatrix}
		n - 1 + \alpha_n
		\\
		n - 1
		\end{pmatrix}\begin{pmatrix}
		n + \alpha\nn
		\\
		n
		\end{pmatrix}
		\\
		= & 0.
		\end{align}
		Here the last equality is from the definition of $\alpha_n$ in \eqref{eq:alphabetax}. 
		When $x\to\infty$, note that the coefficient of the leading term $x^{2n}$ in $\Delta_n(x)$ is 
		\begin{align}
		&\dfrac{(2n+\alpha_n+\beta_n+1)!}{2^nn!(n+\alpha_n+\beta_n+1)!}\dfrac{(2n+\alpha\nn+\beta\nn+1)!}{2^nn!(n+\alpha\nn+\beta\nn+1)!} 
		\\
		& - \dfrac{(2n+\alpha_n+\beta_n-1)!}{2^{n-1}(n-1)!(n+\alpha_n+\beta_n)!}\dfrac{(2n+\alpha\nn+\beta\nn+3)!}{2^{n+1}(n+1)!(n+\alpha\nn+\beta\nn+2)!}
		\\
		= &\dfrac{(2n+\alpha_n+\beta_n-1)!(2n+\alpha\nn+\beta\nn+1)!}{2^{2n}(n-1)!n!(n+\alpha_n+\beta_n)!(n+\alpha\nn+\beta\nn+1)!}\cdot
		\\
		&\ \ \ \left[\dfrac{(2n+\alpha_n+\beta_n)(2n+\alpha_n+\beta_n+1)}{n(n+\alpha_n+\beta_n+1)} - \dfrac{(2n+\alpha\nn+\beta\nn+2)(2n+\alpha\nn+\beta\nn+3)}{(n+1)(n+\alpha\nn+\beta\nn+2)} \right],
		\end{align}
		in which substituting \eqref{eq:alphabetax} we have
		\begin{align}
		& \dfrac{(2n+\alpha_n+\beta_n)(2n+\alpha_n+\beta_n+1)}{n(n+\alpha_n+\beta_n+1)} - \dfrac{(2n+\alpha\nn+\beta\nn+2)(2n+\alpha\nn+\beta\nn+3)}{(n+1)(n+\alpha\nn+\beta\nn+2)}
		\\
		= & \dfrac{(2+a+b)((2+a+b)n+1)}{(1+a+b)n+1} - \dfrac{(2+a+b)((2+a+b)(n+1)+1)}{(1+a+b)(n+1)+1} 
		\\
		= &(2+a+b)\left[\dfrac{1}{1+a+b+\dfrac{1}{n}} - \dfrac{1}{1+a+b+\dfrac{1}{n+1}} \right] <0.
		\end{align}
		Thus we have $\lim_{x\to \infty}\Delta_n(x)=-\infty$. 
		Therefore, we conclude that $f_n^{a,b}(x)< 0$ for all $x> 1$, and hence $\Delta_n(x)< 0$ as well.
	\end{proof}
\end{thm}

\section{Concluding remarks}

In this note we prove a proof of Tu\'ran-like inequality for Jacobi polynomials of form $P_n^{\alpha_n,\beta_n}(x)$. The difference between our result \eqref{eq:Turan} and previous ones in the literature is that in our case $\alpha_n$ and $\beta_n$ are not constants but depends linearly on $n$. While our result applies for all variables $x\ge 1$, it will be interesting to see if similar inequality holds for $x$ in $[-1,1]$ and $(-\infty,-1]$.


\bibliography{yuyuan}

\end{document}